\documentclass{article}
\usepackage[centertags]{amsmath}
\usepackage{amsfonts}
\usepackage{amssymb}
\usepackage{amsthm}
\usepackage{newlfont}
\newcommand{\co}{\mathrm{const}}
\newcommand{\e}{\epsilon}
\newcommand{\la}{\lambda}
\newcommand{\G}{\Gamma}
\begin{document}
\title{Finite extinction time for the solutions to the Ricci flow
on certain three-manifolds}\author{Grisha
Perelman\thanks{St.Petersburg branch of Steklov Mathematical
Institute, Fontanka 27, St.Petersburg 191023, Russia. Email:
perelman@pdmi.ras.ru or perelman@math.sunysb.edu }} \maketitle
\par In our previous paper we constructed complete solutions to
the Ricci flow with surgery for arbitrary initial riemannian
metric on a (closed, oriented) three-manifold [P,6.1], and used
the behavior of such solutions to classify three-manifolds into
three types [P,8.2]. In particular, the first type consisted of
those manifolds, whose prime factors are diffeomorphic copies of
spherical space forms and $\mathbb{S}^2\times\mathbb{S}^1;$ they
were characterized by the property that they admit  metrics, that
give rise to solutions to the Ricci flow with surgery, which
become extinct in finite time. While this classification was
sufficient to answer topological questions, an analytical question
of significant independent interest remained open, namely, whether
the solution becomes extinct in finite time for every initial
metric on a manifold of this type.
\par In this note we prove that this is indeed the case. Our
argument (in conjunction with [P,\S 1-5]) also gives a direct
proof of the so called "elliptization conjecture". It turns out
that it does not require any substantially new ideas: we use only
a version of the least area disk argument from [H,\S 11] and a
regularization of the curve shortening flow from [A-G].
\section{Finite time extinction}
\par {\bf 1.1 Theorem.} {\it Let $M$ be a closed oriented
three-manifold, whose prime decomposition contains no aspherical
factors. Then for any initial metric on $M$ the solution to the
Ricci flow with surgery becomes extinct in finite time.}
\par {\it Proof for irreducible $M$.} Let $\Lambda M$ denote the
space of all contractible loops in $C^1(\mathbb{S}^1\to M).$ Given
a riemannian metric $g$ on $M$ and $c\in\Lambda M,$ define
$A(c,g)$ to be the infimum of the areas of all lipschitz maps from
$\mathbb{D}^2$ to $M,$ whose restriction to
$\partial\mathbb{D}^2=\mathbb{S}^1$ is $c.$ For a family
$\Gamma\subset\Lambda M$ let $A(\Gamma ,g)$ be the supremum of
$A(c,g)$ over all $c\in\Gamma.$ Finally, for a nontrivial homotopy
class $\alpha\in\pi_{\ast}(\Lambda M,M)$ let $A(\alpha ,g)$ be the
infimum of $A(\Gamma ,g)$ over all $\Gamma\in\alpha.$ Since $M$ is
not aspherical, it follows from a classical (and elementary)
result of Serre that such a nontrivial homotopy class exists.
\par {\bf 1.2 Lemma.} (cf. [H,\S 11]) {\it If $g^t$ is a smooth
solution to the Ricci flow, then for any $\alpha$ the rate of
change of the function $A^t=A(\alpha ,g^t)$ satisfies the estimate
$$ \frac{d}{dt}A^t\le-2\pi-\frac{1}{2}R^t_{\mathrm{min}}A^t $$
(in the sense of the $\mathrm{lim \  sup}$ of the forward
difference quotients), where $R^t_{\mathrm{min}}$ denotes the
minimum of the scalar curvature of the metric $g^t.$ }
\par A rigorous proof of this lemma will be given in \S 3, but the
idea is simple and can be explained here. Let us assume that at
time $t$ the value $A^t$ is attained by the family $\Gamma,$ such
that the loops $c\in\Gamma$ where $A(c,g^t)$ is close to $A^t$ are
embedded and sufficiently smooth. For each such $c$ consider the
minimal disk $D_c$ with boundary $c$ and with area $A(c,g^t).$ Now
let the metric evolve by the Ricci flow and let the curves $c$
evolve by the curve shortening flow (which moves every point of
the curve in the direction of its curvature vector at this point)
with the same time parameter. Then the rate of change of the area
of $D_c$ can be computed as
$$\int_{D_c}{(-\mathrm{Tr(Ric^T)})}\ \ + \int_c{(-k_g)}$$
where $\mathrm{Ric^T}$ is the Ricci tensor of $M$ restricted to
the tangent plane of $D_c,$ and $k_g$ is the geodesic curvature of
$c$ with respect to $D_c$ (cf. [A-G, Lemma 3.2]). In three
dimensions the first integrand equals
$-\frac{1}{2}R-(K-\mathrm{det \ II}),$ where $K$ is the intrinsic
curvature of $D_c$ and $\mathrm{det \ II},$ the determinant of the
second fundamental form, is nonpositive, because $D_c$ is minimal.
Thus, the rate of change of the area of $D_c$ can be estimated
from above by $$ \int_{D_c}(-\frac{1}{2}R-K) \ + \int_c (-k_g) =
\int_{D_c}(-\frac{1}{2}R) \ -2\pi $$ by the Gauss-Bonnet theorem,
and the statement of the lemma follows.
\par The problem with this argument is that if $\Gamma$ contains
curves, which are not immersed (for instance, a curve could pass
an arc once in one direction and then make an about turn and pass
the same arc in the opposite direction), then it is not clear how
to define curve shortening flow so that it would be continuous
both in the time parameter and in the family parameter. In \S 3
we'll explain how to circumvent this difficulty, essentially by
adding one dimension to the ambient manifold. This regularization
of the curve shortening flow has been worked out by Altschuler and
Grayson [A-G] (who were interested in approximating the singular
curve shortening flow on the plane and obtained for that case more
precise results than what we need).
\par {\bf 1.3} Now consider the solution to the Ricci flow with
surgery. Since $M$ is assumed irreducible, the surgeries are
topologically trivial, that is one of the components of the
post-surgery manifold is diffeomorphic to the pre-surgery
manifold, and all the others are spheres. Moreover, by the
construction of the surgery [P,4.4], the diffeomorphism from the
pre-surgery manifold to the post-surgery one can be chosen to be
distance non-increasing ( more precisely, $(1+\xi)$-lipschitz,
where $\xi>0$ can be made as small as we like). It follows that
the conclusion of the lemma above holds for the solutions to the
Ricci flow with surgery as well.
\par Now recall that the evolution equation for the scalar
curvature $$ \frac{d}{dt}R=\triangle R +
2|\mathrm{Ric}|^2=\triangle
R+\frac{2}{3}R^2+2|\mathrm{Ric}^{\circ}|^2 $$ implies the estimate
$R^t_{\mathrm{min}}\ge -\frac{3}{2}\frac{1}{t+\co}.$ It follows
that $\hat{A}^t=\frac{A^t}{t+\co}$ satisfies
$\frac{d}{dt}\hat{A}^t\le -\frac{2\pi}{t+\co},$ which implies
finite extinction time since the right hand side is non-integrable
at infinity whereas $\hat{A}^t$ can not become negative.
\par {\bf 1.4} {\it Remark.} The finite time extinction result for
irreducible non-aspherical manifolds already implies (in
conjuction with the work in [P,\S 1-5] and the Kneser finiteness
theorem) the so called "elliptization conjecture", claiming that a
closed manifold with finite fundamental group is diffeomorphic to
a spherical space form. The analysis of the long time behavior in
[P,\S 6-8] is not needed in this case; moreover the argument in
[P,\S 5] can be slightly simplified, replacing the sequences $r_j,
\kappa_j, \bar{\delta}_j$ by single values $r, \kappa,
\bar{\delta},$ since we already have an upper bound on the
extinction time in terms of the initial metric.
\par In fact, we can even avoid the use of the Kneser theorem.
Indeed, if we start from an initial metric on a homotopy sphere
(not assumed irreducible), then at each surgery time we have
(almost) distance non-increasing homotopy equivalences from the
pre-surgery manifold to each of the post-surgery components, and
this is enough to keep track of the nontrivial relative homotopy
class of the loop space.
\par {\bf 1.5} {\it Proof of theorem 1.1 for general $M$.}
The Kneser theorem implies that our solution undergoes only
finitely many topologically nontrivial surgeries, so from some
time $T$ on all the surgeries are trivial. Moreover, by the Milnor
uniqueness theorem, each component at time $T$ satisfies the
assumption of the theorem. Since we already know from 1.4 that
there can not be any simply connected prime factors, it follows
that every such component is either irreducible, or has nontrivial
$\pi_2;$ in either case the proof in 1.1-1.3 works.
\section{Preliminaries on the curve shortening flow}
\par In this section we rather closely follow [A-G].
\par {\bf 2.1} Let $M$ be a closed $n$-dimensional manifold, $n\ge
3,$ and let $g^t$ be a smooth family of riemannian metrics on $M$
evolving by the Ricci flow on a finite time interval $[t_0,t_1].$
It is known [B] that $g^t$ for $t>t_0$ are real analytic. Let
$c^t$ be a solution to the curve shortening flow in $(M,g^t),$
that is $c^t$ satisfies the equation $\frac{d}{dt}c^t(x)=H^t(x),$
where $x$ is the parameter on $\mathbb{S}^1,$ and $H^t$ is the
curvature vector field of $c^t$ with respect to $g^t.$ It is known
[G-H] that for any smoothly immersed initial curve $c$ the
solution $c^t$ exists on some time interval $[t_0,t_1'),$ each
$c^t$ for $t>t_0$ is an analytic immersed curve, and either
$t_1'=t_1,$ or the curvature $k^t=g^t(H^t,H^t)^{\frac{1}{2}}$ is
unbounded when $t\to t_1'.$
\par Denote by $X^t$ the tangent vector field to $c^t,$ and let
$S^t=g^t(X^t,X^t)^{-\frac{1}{2}}X^t$ be the unit tangent vector
field; then $H=\nabla_S S$ (from now on we drop the superscript
$t$ except where this omission can cause confusion). We compute
\begin{equation}
\frac{d}{dt}g(X,X)=-2\mathrm{Ric}(X,X)-2g(X,X)k^2, \end{equation}
which implies \begin{equation} [H,S]=(k^2+\mathrm{Ric}(S,S))S
\end{equation} Now we can compute \begin{equation}
\frac{d}{dt}k^2=(k^2)''-2g((\nabla_S H)^{\perp},(\nabla_S
H)^{\perp})+2k^4 + ..., \end{equation} where primes denote
differentiation with respect to the arclength parameter $s,$ and
where dots stand for the terms containing the curvature tensor of
$g,$ which can be estimated in absolute value by
$\co\cdot(k^2+k).$ Thus the curvature $k$ satisfies
\begin{equation} \frac{d}{dt}k\le k''+k^3+\co\cdot(k+1)
\end{equation} Now it follows from (1) and (4) that the length $L$
and the total curvature $\Theta=\int kds$ satisfy \begin{equation}
\frac{d}{dt}L\le \int (\co \ -k^2)ds, \end{equation}
\begin{equation} \frac{d}{dt}\Theta\le \int \co\cdot (k+1)ds
\end{equation}
In particular, both quantities can grow at most exponentially in
$t$ (they would be non-increasing in a flat manifold).
\par {\bf 2.2} In general the curvature of $c^t$ may concentrate
near certain points, creating singularities. However, if we know
that this does not happen at some time $t^{\ast},$ then we can
estimate the curvature and higher derivatives at times shortly
thereafter. More precisely, there exist constants $\e, C_1,
C_2,...$ (which may depend on the curvatures of the ambient space
and their derivatives, but are independent of $c^t$), such that if
at time $t^{\ast}$ for some $r>0$ the length of $c^t$ is at least
$r$ and the total curvature of each arc of length $r$ does not
exceed $\e,$ then for every $t\in (t^{\ast},t^{\ast}+\e r^2)$ the
curvature $k$ and higher derivatives satisfy the estimates
$k^2=g(H,H)\le C_0 (t-t^{\ast})^{-1},\ \  g(\nabla_S H,\nabla_S
H)\le C_1 (t-t^{\ast})^{-2},...\ \ $ This can be proved by
adapting the arguments of Ecker and Huisken [E-Hu]; see also
[A-G,\S 4].
\par {\bf 2.3} Now suppose that our manifold $(M,g^t)$ is a metric
product $(\bar{M},\bar{g}^t)\times \mathbb{S}^1_{\la},$ where the
second factor is the circle of constant length $\la;$ let $U$
denote the unit tangent vector field to this factor. Then
$u=g(S,U)$ satisfies the evolution equation \begin{equation}
\frac{d}{dt}u=u''+(k^2+\mathrm{Ric}(S,S))u
\end{equation} \par Assume that $u$ was strictly
positive everywhere at time $t_0$ (in this case the curve is
called a ramp). Then it will remain positive and bounded away from
zero as long as the solution exists. Now combining (4) and (7) we
can estimate the right hand side of the evolution equation for the
ratio $\frac{k}{u}$ and conclude that this ratio, and hence the
curvature $k,$ stays bounded (see [A-G,\S 2]). It follows that
$c^t$ is defined on the whole interval $[t_0,t_1].$
\par {\bf 2.4} Assume now that we have two ramp solutions $c_1^t,
c_2^t,$ each winding once around the $\mathbb{S}^1_{\la}$ factor.
Let $\mu^t$ be the infimum of the areas of the annuli with
boundary $c_1^t\cup c_2^t.$ Then \begin{equation}
\frac{d}{dt}\mu^t\le (2n-1)|\mathrm{Rm}^t|\mu^t, \end{equation}
where $|\mathrm{Rm}^t|$ denotes a bound on the absolute value of
sectional curvatures of $g^t.$ Indeed, the curves $c^t_1$ and
$c^t_2,$ being ramps, are embedded and without substantial loss of
generality we may assume them to be disjoint. In this case the
results of Morrey [M] and Hildebrandt [Hi] yield an analytic
minimal annulus $A,$ immersed, except at most finitely many branch
points, with prescribed boundary and with area $\mu.$ The rate of
change of the area of $A$ can be computed as $$ \int_A
(-\mathrm{Tr}(\mathrm{Ric}^T)) +\int_{\partial A} (-k_g) \le
\int_A (-\mathrm{Tr}(\mathrm{Ric}^T)+K)$$ $$ \le \int_A
(-\mathrm{Tr}(\mathrm{Ric}^T)+\mathrm{Rm}^T) \le
(2n-1)|\mathrm{Rm}|\mu ,$$ where the first inequality comes from
the Gauss-Bonnet theorem, with possible contribution of the branch
points, and the second one is due to the fact that a minimal
surface has nonpositive extrinsic curvature with respect to any
normal vector.
\par {\bf 2.5} The estimate (8) implies that $\mu^t$ can grow at
most exponentially; in particular, if $c^t_1$ and $ c^t_2$ were
very close at time $t_0,$ then they would be close for all $t\in
[t_0,t_1]$ in the sense of minimal annulus area. In general this
does not imply that the lengths of the curves are also close.
However, an elementary argument shows that if $\e>0$ is small
then, given any $r>0, $ one can find $\bar{\mu},$ depending only
on $r$ and on upper bound for sectional curvatures of the ambient
space, such that if the length of $c^t_1$ is at least $r,$ each
arc of $c^t_1$ with length $r$ has total curvature at most $\e,$
and $\mu^t\le\bar{\mu},$ then $L(c^t_2)\ge (1-100\e)L(c^t_1).$
\section{Proof of lemma 1.2}
\par  {\bf 3.1} In this section we prove the following statement
\par {\it Let $M$ be a closed three-manifold, and let $(M,g^t)$ be a
smooth solution to the Ricci flow on a finite time interval
$[t_0,t_1].$ Suppose that $\G\subset \Lambda M$ is a compact
family. Then for any $\xi>0$ one can construct a continuous
deformation $\G^t, t\in [t_0,t_1], \G^{t_0}=\G,$ such that for
each curve $c\in\G$  either the value $A(c^{t_1},g^{t_1})$ is
bounded from above by $\xi$ plus the value at $t=t_1$ of the
solution to the ODE
$\frac{d}{dt}w(t)=-2\pi-\frac{1}{2}R^t_{\mathrm{min}}w(t)$ with
the initial data $\ \ w(t_0)=A(c^{t_0},g^{t_0}),\ \ $ or
$L(c^{t_1})\le \xi;$ moreover, if $c$ was a constant map, then all
$c^t$ are constant maps.}
\par It is clear that our statement implies lemma 1.2, because a
family consisting of very short loops can not represent a
nontrivial relative homotopy class.
\par {\bf 3.2} As a first step of the proof of the statement we
can  replace $\G$ by a family, which consists of piecewise
geodesic loops with some large fixed number of vertices and with
each segment reparametrized in some standard way to make the
parametrizations of the whole curves twice continuously
differentiable.
\par Now consider the manifold $M_{\la}=M\times\mathbb{S}^1_{\la},
0<\la<1,$ and for each $c\in\G$ consider the smooth embedded
closed curve $c_{\la}$ such that $p_1c_{\la}(x)=c(x)$ and
$p_2c_{\la}(x)=\la x\  \mathrm{mod}\  \la,$ where $p_1$ and $p_2$
are projections of $M_{\la}$ to the first and second factor
respectively, and $x$ is the parameter of the curve $c$ on the
standard circle of length one. Using 2.3 we can construct a
solution $c_{\la}^t, t\in [t_0,t_1]$ to the curve shortening flow
with initial data $c_{\la}.$ The required deformation will be
obtained as $\G^t=p_1\G^t_{\la}$ (where $\G^t_{\la}$ denotes the
family consisting of $c^t_{\la}$) for certain sufficiently small
$\la>0.$ We'll verify that an appropriate $\la$ can be found for
each individual curve $c,$ or for any finite number of them, and
then show that if our $\la$ works for all elements of a $\mu$-net
in $\G,$ for sufficiently small $\mu>0,$ then it works for all
elements of $\G.$
\par {\bf 3.3} In the following estimates we shall denote by $C$
large constants that may depend on metrics $g^t,$ family $\G$ and
$\xi,$ but are independent of $\la ,\mu $ and a particular curve
$c.$
\par The first step in 3.2 implies that the lengths and total
curvatures of $c_{\la}$ are uniformly bounded, so by 2.1 the same
is true for all $c^t_{\la}.$ It follows that the area swept by
$c^t_{\la}, t\in [t',t'']\subset [t_0,t_1]$ is bounded above by
$C(t''-t'),$ and therefore we have the estimates
$A(p_1c^t_{\la},g^t)\le C,
A(p_1c^{t''}_{\la},g^{t''})-A(p_1c^{t'}_{\la},g^{t'})\le
C(t''-t').$
\par {\bf 3.4} It follows from (5) that $\int_{t_0}^{t_1} \int k^2
dsdt \le C$ for any $c^t_{\la}.$ Fix some large constant $B,$ to
be chosen later. Then there is a subset $I_B(c_{\la})\subset
[t_0,t_1]$ of measure at least $t_1-t_0-CB^{-1}$ where $\int k^2
ds\le B,$ hence $\int kds \le \e$ on any arc of length $\le \e^2
B^{-1}.$ Assuming that $c^t_{\la}$ are at least that long, we can
apply 2.2 and construct another subset $J_B(c_{\la})\subset
[t_0,t_1]$ of measure at least $t_1-t_0-CB^{-1},$ consisting of
finitely many intervals of measure at least $C^{-1}B^{-2}$ each,
such that for any $t\in J_B(c_{\la})$ we have pointwise estimates
on $c^t_{\la}$ for curvature and higher derivatives, of the form
$k\le CB,...$
\par Now fix $c, B,$ and consider any sequence of $\la\to 0.$
Assume again that the lengths of $c^t_{\la}$ are bounded below by
$\e^2 B^{-1},$ at least for $t\in [t_0,t_2],$ where
$t_2=t_1-B^{-1}.$ Then an elementary argument shows that we can
find a subsequence $\Lambda_c$ and a subset $J_B(c)\subset
[t_0,t_2]$ of measure at least $t_1-t_0-CB^{-1},$ consisting of
finitely many intervals, such that $J_B(c)\subset J_B(c_{\la})$
for all $\la\in \Lambda_c.$ It follows that on every interval of
$J_B(c)$ the curve shortening flows $c^t_{\la}$ smoothly converge
(as $\la\to 0$ in some subsequence of $\Lambda_c$ )  to a curve
shortening flow in $M.$
\par Let $w_c(t)$ be the solution of the ODE
$\frac{d}{dt}w_c(t)=-2\pi-\frac{1}{2}R^t_{\mathrm{min}}w_c(t)$
with initial data $w_c(t_0)=A(c,g^{t_0}).$ Then for sufficiently
small $\la\in\Lambda_c$ we have $A(p_1c^t_{\la},g^t)\le
w_c(t)+\frac{1}{2}\xi$ provided that $B>C\xi^{-1}.$ Indeed, on the
intervals of $J_B(c)$ we can estimate the change of $A$ for the
limit flow using the minimal disk argument as in 1.2, and this
implies the corresponding estimate for $p_1c^t_{\la}$ if
$\la\in\Lambda_c$ is small enough, whereas for the intervals of
the complement of $J_B(c)$ we can use the estimate in 3.3.
\par On the other hand, if our assumption on the lower bound for
lengths does not hold, then it follows from (5) that
$L(c^{t_2}_{\la})\le CB^{-1}\le\frac{1}{2}\xi.$
\par {\bf 3.5} Now apply   the previous argument to all elements
of some finite $\mu$-net $\hat{\G}\subset \G$ for small $\mu>0$ to
be determined later. We get a $\la>0$ such that for each
$\hat{c}\in\hat{\G}$ either $A(p_1\hat{c}_{\la}^{t_1},g^{t_1})\le
w_{\hat{c}}(t_1)+\frac{1}{2}\xi$ or $L(\hat{c}_{\la}^{t_2})\le
\frac{1}{2}\xi.$ Now for any curve $c\in\G$ pick a curve
$\hat{c}\in\hat{\G},$ $\mu$-close to $c,$ and apply the result of
2.4. It follows that if $A(p_1\hat{c}_{\la}^{t_1},g^{t_1})\le
w_{\hat{c}}(t_1)+\frac{1}{2}\xi$ and $\mu \le C^{-1}\xi,$ then
$A(p_1c_{\la}^{t_1},g^{t_1})\le w_c(t_1)+\xi.$ On the other hand,
if $L(\hat{c}_{\la}^{t_2})\le \frac{1}{2}\xi,$ then we can
conclude that $L(c_{\la}^{t_1})\le \xi$ provided that $\mu>0$ is
small enough in comparison with $\xi$ and $B^{-1}.$ Indeed, if
$L(c_{\la}^{t_1})>\xi,$ then $L(c_{\la}^t)>\frac{3}{4}\xi$ for all
$t\in [t_2,t_1];$ on the other hand, using (5) we can find a $t\in
[t_2,t_1],$ such that $\int k^2ds \le CB$ for $c^t_{\la};$ hence,
applying 2.5, we get $L(\hat{c}_{\la}^t)>\frac{2}{3}\xi$ for this
$t,$ which is incompatible with $L(\hat{c}_{\la}^{t_2})\le
\frac{1}{2}\xi.$ The proof of the statement 3.1 is complete.
\section*{References}
\ \ \ [A-G] S.Altschuler, M.Grayson Shortening space curves and
flow through singularities. Jour. Diff. Geom. 35 (1992), 283-298.
\par [B] S.Bando Real analyticity of solutions of Hamilton's
equation. Math. Zeit. 195 (1987), 93-97.
\par [E-Hu] K.Ecker, G.Huisken Interior estimates for
hypersurfaces moving by mean curvature. Invent. Math. 105 (1991),
547-569.
\par [G-H] M.Gage, R.S.Hamilton The heat equation shrinking convex
plane curves. Jour. Diff. Geom. 23 (1986), 69-96.
\par [H] R.S.Hamilton Non-singular solutions of the Ricci flow on
three-manifolds. Commun. Anal. Geom. 7 (1999), 695-729.
\par [Hi] S.Hildebrandt Boundary behavior of minimal surfaces.
Arch. Rat. Mech. Anal. 35 (1969), 47-82.
\par [M] C.B.Morrey The problem of Plateau on a riemannian
manifold. Ann. Math. 49 (1948), 807-851.
\par [P] G.Perelman Ricci flow with surgery on three-manifolds. \par arXiv:math.DG/0303109 v1

\end{document}